\documentclass[11pt]{amsart}
\hoffset= - 0.75 in

\usepackage{amsfonts, amsmath, amssymb, amsthm}
\usepackage{latexsym}
\usepackage{graphics}
\newtheorem{thm}{Theorem}

\newtheorem{ex}[thm]{Example}

\newtheorem{lem}[thm]{Lemma}
\newtheorem{cor}[thm]{Corollary}
\newtheorem{prop}[thm]{Proposition}

\newcommand{\RR}{\mathbb{R}}

\newcommand{\TP}{\mathbb{TP}}
\newcommand{\bin}[2]{{#1\choose #2}}

\begin{document}

\title{Tropical Convexity}

\author[Develin and Sturmfels]{Mike Develin and Bernd Sturmfels}

\address{Mike Develin, American Institute of Mathematics, 360 Portage Ave., Palo Alto, CA 94306-2244, USA}
\email{develin@post.harvard.edu}
\address{Bernd Sturmfels, Department of Mathematics, University
of California, Berkeley, CA 94720, USA}
\email{bernd@math.berkeley.edu}

\begin{abstract}
The notions of convexity and convex polytopes are 
introduced in the setting of tropical geometry.
Combinatorial types of tropical polytopes
are shown to be in bijection with regular
triangulations of products of two simplices.
Applications to phylogenetic trees are discussed.
\end{abstract}

\maketitle

\section{Introduction}

The \emph{tropical semiring} $(\RR, \oplus, \odot)$ is the set
of real numbers with the arithmetic operations of
\emph{tropical addition}, which is taking the minimum of two numbers,
and \emph{tropical multiplication}, which is ordinary addition.
Thus the two arithmetic operations are defined as follows:
$$ a \oplus b \,\, := \,\, {\rm min} (a,b) \qquad
\hbox{and} \qquad a \odot b \,\, := \,\, a + b . $$
The $n$-dimensional space $\RR^n$  is a semimodule over the
tropical semiring, with tropical addition
$$  (x_1,\ldots,x_n) \oplus (y_1,\ldots,y_n) \quad = \quad
(x_1 \oplus y_1, \ldots, x_n \oplus y_n) , $$
and tropical scalar multiplication
$$ c \odot (x_1,x_2,\ldots,x_n) \quad = \quad (c \odot x_1,
c \odot x_2, \ldots, c \odot x_n) . $$
The semiring $(\RR, \oplus, \odot)$ and its semimodule
$\RR^n$ obey the usual distributive and associative laws.

The purpose of this paper is to propose a tropical theory
of convex polytopes.
Convexity in arbitrary idempotent semimodules was introduced
by Cohen, Gaubert  and Quadrat \cite{CGQ}
and  Litvinov, Maslov and Shpiz \cite{LMS}. 
Some of our results (such as Theorem \ref{polar}
and Propositions \ref{intersectpolyt} and
\ref{vertexunique}) are known in a different
guise in idempotent analysis.
Our objective is to  provide
a combinatorial approach to convexity 
in the tropical semiring which is
consistent with the recent developments in
tropical algebraic geometry (see \cite{Mi}, \cite{RGST}, \cite{SS}).
The connection to tropical methods in
representation theory  (see \cite{K},  \cite{NY}) is
less clear and deserves further study.

There are many notions of discrete convexity in
the computational geometry literature, but none of them
seems to be quite like tropical convexity. For instance,
the notion of \emph{directional convexity} studied by
Matou\v sek \cite{Mat} has similar features
but it is different and much harder to compute with.

A subset  $S$ of $ \RR^n$ is called
\emph{tropically convex} if
the set $S$ contains the point
$\,a \odot x \,\oplus \, b \odot y \,$
 for all $x,y\in S$ and all $ \, a,b\in \RR $.
The \emph{tropical convex hull} of a given subset $V \subset \RR^n$ is the
smallest tropically convex subset of $\RR^n$ which contains $V$.
We shall see in Proposition \ref{tstc}
that the tropical convex hull of $V$
coincides with the
set of all tropical linear combinations
\begin{equation}
\label{alltropcomb}
  a_1 \odot v_1 \,\oplus \,
a_2 \odot v_2 \,\oplus \, \cdots \,\oplus \, a_r \odot v_r \, ,
\qquad \hbox{where} \,\,\, v_1,\ldots,v_r \in V \,\,\,
\hbox{and} \,\,\, a_1,\ldots,a_r \in \RR.
\end{equation}
Any tropically convex subset $S$ of  $\RR^n$ is
closed under tropical scalar multiplication,
$\, \RR \odot S \, \subseteq \,S$. In other words,
if $\, x \in S \,$ then
$\, x + \lambda (1,\ldots,1) \, \in \,S \,$
for all $\lambda \in \RR$.
We will therefore identify the tropically convex set $S$
with its image in the
$(n-1)$-dimensional  \emph{tropical projective space}
$$\TP^{n-1} \quad = \quad \RR^n/(1,\ldots,1)\RR . $$
Basic properties of (tropically) convex subsets in $\TP^{n-1}$
will be presented in Section 2.
In Section 3 we introduce tropical polytopes
and study their combinatorial structure.
 A \emph{tropical polytope} is the tropical convex hull
of a finite subset $V$ in $\TP^{n-1}$.
Every tropical polytope is a finite union of
convex polytopes in the usual sense: given a set $V=\{v_1,\ldots,v_n\}$, their convex hull has a natural 
decomposition as a polyhedral complex, which we call the \emph{tropical complex} generated by $V$.
The following  main result will be proved in Section 4:

\begin{thm} \label{main}
The combinatorial types of tropical complexes generated by a set of 
$r$ vertices in $\TP^{n-1}$ are in natural bijection with the
regular polyhedral subdivisions of the
product of two simplices $\Delta_{n-1} \times \Delta_{r-1}$.
\end{thm}

This implies a remarkable duality between
tropical $(n-1)$-polytopes with $r$ vertices and
tropical $(r-1)$-polytopes with $n$ vertices.  Another consequence
of Theorem~\ref{main} is a formula for the $f$-vector of a generic tropical complex.
In Section 5 we discuss applications of
tropical convexity to phylogenetic analysis,
extending known results on injective hulls of 
finite metric spaces
(cf.~\cite{DHM}, \cite{DMT}, \cite{DT} and \cite{SS}).

\section{Tropically convex sets}

\begin{figure}
 \input{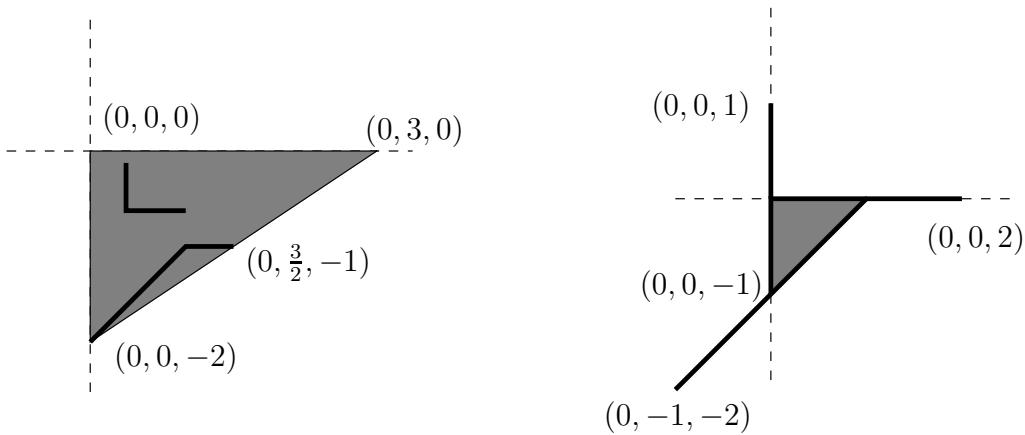}
\caption{\label{tropconv}
Tropical convex sets and tropical line segments in $\TP^2$.
}
\end{figure}

We begin with two pictures of tropical convex sets in
the tropical plane $\TP^2$. A point $(x_1,x_2,x_3)
\in \TP^2$ is represented by drawing
the point with coordinates $\,(x_2-x_1,x_3-x_1)\,$
in the plane of the paper. The triangle on the left hand side in
Figure~\ref{tropconv}
is tropically convex, but it is not a tropical polytope
because it is not the tropical convex hull
of finitely many points. The thick edges indicate
two tropical line segments. The picture on the
right hand side is a \emph{tropical triangle}, namely,
it is the tropical convex hull of the three points
$\, (0,0,1), \, (0,2,0) $ and $(0,-1,-2)$ in the tropical plane
$\TP^2$. The thick
edges represent the tropical segments connecting
any two of these three points.

We next show that
tropical convex sets enjoy many of the features of ordinary convex sets.

\begin{thm} \label{LotsOfStuff}
The intersection of
two tropically convex sets  in $\RR^n$ or in $\TP^{n-1}$
is tropically convex.
The projection of a tropically convex set onto a coordinate hyperplane is tropically convex.
The ordinary hyperplane $\{x_i-x_j=l\}$ is tropically convex, and the
projection map from this hyperplane to $\RR^{n-1}$ given by eliminating $x_i$ is an isomorphism of tropical semimodules.
Tropically convex sets are contractible
spaces. The Cartesian product of two tropically convex sets is tropically convex.
\end{thm}

\begin{proof}
We prove the statements in the order given. If $S$ 
and $T$ are tropically convex, then for any two points $x,y\in S\cap T$, both
$S$ and $T$ contain the tropical
line segment between $x$ and $y$, and consequently so does $S\cap T$. Therefore $S\cap T$ is
tropically convex by definition.

Suppose $S$ is a tropically convex set   in $\RR^n$.
 We wish to show that the image of $S$ under the  coordinate projection
$\,\phi : \RR^n \rightarrow \RR^{n-1}, \,
(x_1,x_2,\ldots,x_n) \mapsto (x_2,\ldots,x_n)\,$ is a tropically convex
subset of $\RR^{n-1}$. If $x,y \in S$ then we have the obvious identity
$$ \phi \bigl(c \odot x \,\oplus \,d \odot y \bigr)
\quad = \quad c \odot \phi(x) \, \oplus \,d \odot \phi(y). $$
This means that $\,\phi\,$
is a homomorphism of tropical semimodules.
 Therefore, if $S$ contains the tropical line segment between $x$ and
$y$, then $\phi(S)$ contains the tropical
line segment between $\phi(x)$ and $\phi(y)$ and hence is
tropically convex.  The same holds for the induced map
$\,\phi : \TP^{n-1} \rightarrow \TP^{n-2}$.

Most ordinary hyperplanes in $\RR^n$ are not tropically convex, but
we are claiming that hyperplanes of the special form
$\,x_i-x_j=k \,$ are tropically convex. If $x$ and $y$ lie in that
hyperplane then $\,x_i-y_i = x_j  - y_j $. This last equation implies
the following identity for any real numbers $c,d \in \RR$:
$$
(c \odot x \,\oplus \, d \odot y)_i \,-\,
(c \odot x \,\oplus \, d \odot y)_j  \quad = \quad
{\rm min}(x_i + c, y_i + d ) \, -\,
{\rm min}(x_j + c, y_j + d ) \quad = \quad k . $$
Hence the tropical line segment between $x$ and $y$ also lies
in the hyperplane $\{x_i  - x_j = k \}$.

Consider the map from $\{x_i  - x_j = k \}$
to $\RR^{n-1}$ given by deleting the $i$-th coordinate.
This map is injective: if two points
differ in the $x_i$ coordinate they must also differ in the
$x_j$ coordinate.  It is clearly surjective because we
can recover an $i$-th coordinate by setting
$\,x_i = x_j + k$. Hence this map is an isomorphism
of $\RR$-vector spaces and it is also an isomorphism
of $(\RR, \oplus, \odot)$-semimodules.

Let $S$ be a tropically convex set in $\RR^n$ or $\TP^{n-1}$.
Consider the family of hyperplanes $H_l = \{x_1-x_2=l\}$ for $l\in \RR$.
We know that the intersection $S\cap H_l$ is tropically convex,
and isomorphic to its (convex) image under the map deleting the first
coordinate. This image is contractible by induction on the dimension $n$ of the ambient space. Therefore, $S\cap H_l$
is contractible. The result then follows from the topological result that if $S$ is connected, which all tropically
convex sets obviously are, and if $S\cap H_l$ is contractible for each $l$, then $S$ itself is also contractible.

Suppose that $S\subset \RR^n$ and $T\subset \RR^m$ are tropically convex.
Our last assertion states that $S\times T$ is a tropically
convex subset of $\RR^{n+m}$. Take any $(x,y)$ and $(x',y')$ in $S\times T$
and $c,d \in \RR$. Then
$$ c \odot (x,y) \,\oplus \,
d \odot (x',y') \quad = \quad
\bigl( \,c \odot x \oplus d \odot x' \, , \,
         c \odot y \oplus d \odot y' \, \bigr) $$
lies in $S \times T$ since $S$ and $T$ are tropically convex.
\end{proof}

We next give a more precise description of
what tropical line segments look like.

\begin{prop}\label{slopes}
The tropical line segment between two points $x$ and $y$ in $\TP^{n-1}$
is the concatenation of at most $n-1$ ordinary line  segments.
The slope of each line segment is a zero-one vector.
\end{prop}

\begin{proof}
After relabeling the coordinates of
 $x = (x_1,,\ldots,x_n) $ and
$y = (y_1,\ldots,y_n)$, we may assume
\begin{equation}
\label{sortthedifferences}
 y_1 - x_1
\,\leq\, y_2 - x_2\, \leq\, \cdots\, \leq \,y_n - x_n .
\end{equation}
The following points lie in the given order
 on the tropical segment between $x$ and $y$:
\begin{eqnarray*}
x \,\,\, = \quad (y_1-x_1) \odot x \,\oplus \,  y \quad = &
\bigl( y_1,\,
 y_1\!-\!x_1 \!+\! x_2,
 y_1\!-\!x_1 \!+\! x_3, \ldots,
 y_1\!- \!x_1 \!+\! x_{n-1},
 y_1\!-\!x_1 \!+\! x_n \bigr) \\
 (y_2-x_2) \odot x \,\oplus \,  y  \quad = &
 \bigl( y_1,\,  y_2, \,
 y_2\!-\!x_2 \!+\! x_3, \ldots,
 y_2\!- \!x_2 \!+\! x_{n-1},
  y_2\!-\!x_2 \!+\! x_n \bigr) \\
 (y_3-x_3) \odot x \,\oplus \,  y  \quad = &
 \bigl( y_1,\,  y_2, \, y_3, \, \ldots,
 y_3\!- \!x_3 \!+\! x_{n-1},
  y_3\!-\!x_3 \!+\! x_n \bigr) \\
 \cdots \cdots \qquad \qquad & \qquad \cdots \cdots  \qquad \cdots \cdots
\\
 (y_{n-1}-x_{n-1}) \odot x \,\oplus \,  y \quad = &
 \bigl( y_1,\,  y_2, \, y_3, \, \ldots, y_{n-1},
 y_{n-1}\!-\!x_{n-1} \!+\! x_n \bigr)  \\
 y \,\,\, = \quad  (y_{n}-x_{n}) \odot x \,\oplus \,  y \quad = &
 \bigl( y_1,\,  y_2, \, y_3, \, \ldots,\, y_{n-1}, \,y_n \bigr).
\end{eqnarray*}
Between any two consecutive points, the tropical
line segment agrees with the ordinary line segment, which has slope
$\,(0,0,\ldots,0,\,1,1,\ldots,1)$.
Hence the tropical line segment between $x$ and $y$
is the concatenation of at most $n-1$ ordinary line segments,
one for each strict inequality in
(\ref{sortthedifferences}).
\end{proof}

This description of tropical segments shows an important
feature of tropical polytopes: their edges
use a limited set of directions. The following
result characterizes the \emph{tropical convex hull}.

\begin{prop}\label{tstc}
The smallest tropically convex subset of $\TP^{n-1}$
which contains a given set $V$ coincides with
the set of all tropical linear combinations (\ref{alltropcomb}).
We denote this set by $\, {\rm tconv}(V)$.
\end{prop}

\begin{proof}
Let $\,x\,= \, \bigoplus_{i=1}^r a_i \odot v_i \,$ be
the point in (\ref{alltropcomb}). If $r \leq 2$ then
$x$ is clearly in the tropical convex hull of $V$.
If $r > 2$ then we write
$\, x \, = \, a_1 \odot v_1 \,\oplus \,
 ( \bigoplus_{i=2}^r a_i \odot v_i)$. The parenthesized vector
lies  the tropical convex hull, by induction on $r$,
and hence so does $x$.
For the converse, consider any two tropical linear combinations
$\,x\,= \, \bigoplus_{i=1}^r c_i \odot v_i \,$ and
$\,y\,= \, \bigoplus_{j=1}^r d_i \odot v_i $.
By the distributive law, $\, a \odot x \,\oplus \,b \odot y \,$ is
also a tropical linear combination of $v_1,\ldots,v_r \in V$.
Hence the set of all tropical linear combinations of $V$
is tropically convex, so it contains the
tropical convex hull of $V$.
\end{proof}

\begin{figure}
 \input{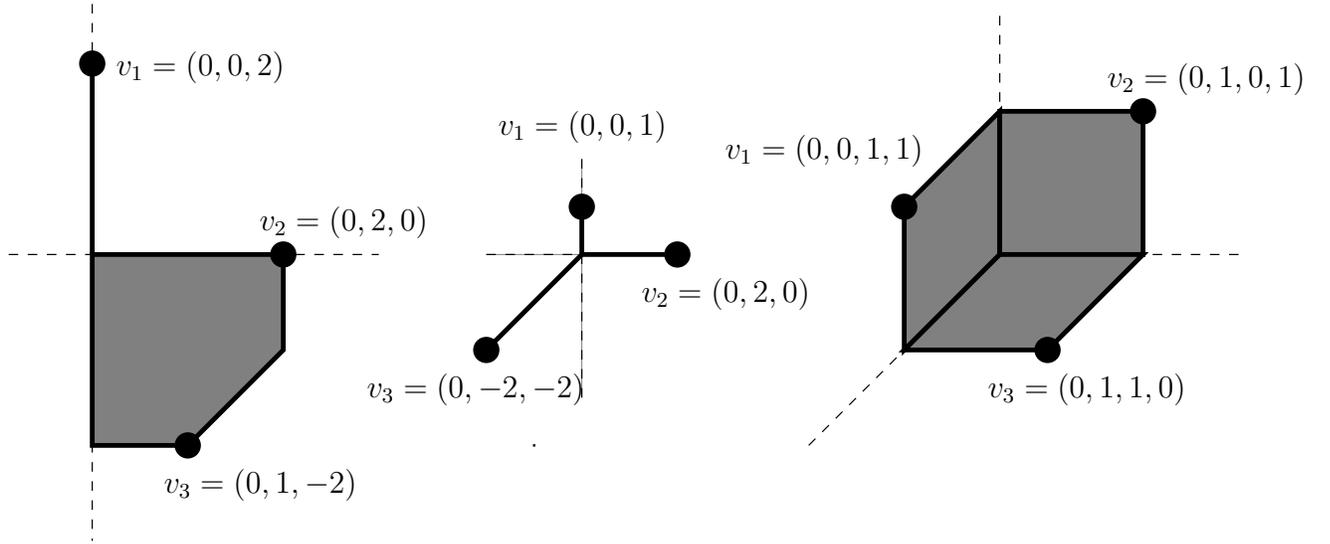}
\caption{\label{throughex}
Three tropical polytopes. The first two live in $\TP^{2}$, the last in $\TP^3$.}
\end{figure}

If $V$ is a finite subset of $\TP^{n-1}$ then $\, {\rm tconv}(V)\,$
is a  \emph{tropical polytope}. In Figure~\ref{throughex}
we see three small examples of tropical polytopes.
The first and second are tropical convex hulls
of three points in $\TP^2$. The third tropical polytope
lies in $\TP^3$ and is  the union of three squares.

One of the basic results in the usual theory of convex polytopes
is Carath\'{e}odory's theorem. This theorem holds in the tropical setting.

\begin{prop}[Tropical Carath\'{e}odory's Theorem]\label{carat}
If $x$ is in the tropical convex hull of a set
of $r$ points $v_i$ in $\TP^{n-1}$, then $x$ is in the tropical
convex hull of at most $n$ of them.
\end{prop}

\begin{proof}
Let $\,x\,= \, \bigoplus_{i=1}^r a_i \odot v_i \,$
and suppose $r > n$. For each coordinate $j \in
\{1,\ldots,n\}$, there exists an index $i \in \{1,\ldots,r \}$ such that
$\,x_j = c_i + v_{ij}$. Take a
subset $I$ of $\{1,\ldots,r\}$ composed of one such $i$ for each $j$.
Then we also have $\,x\,= \, \bigoplus_{i \in I} a_i \odot v_i $,
where $\# (I) \leq n$.
\end{proof}

The basic theory of tropical linear subspaces in $\TP^{n-1}$ was developed in
\cite{RGST} and \cite{SS}. Recall that the \emph{tropical hyperplane}
defined by a tropical linear form
$\,a_1 \odot x_1 \,\oplus \,
a_2 \odot x_2 \,\oplus \,\cdots\,\oplus\,
a_n \odot x_n \,$
consists of all points $x = (x_1,x_2,\ldots,x_n)$ in $\TP^{n-1}$ such that
the following holds (in ordinary arithmetic):
\begin{equation}
\label{hyperplane}
 a_i + x_i \, = \,
a_j + x_j \, = \, {\rm min} \{a_k + x_k \, : \, k =1,\ldots,n \} \quad
\qquad \hbox{for some indices $i \not= j$}.
\end{equation}
Just like in ordinary geometry, hyperplanes are convex sets:

\begin{prop} \label{hypconv}
Tropical hyperplanes in $\TP^{n-1}$ are tropically convex.
\end{prop}

\begin{proof}
Let $H$ be the hyperplane defined by (\ref{hyperplane}).
Suppose that $x$ and $y$ lie in $H$ and consider any
tropical linear combination $\, z \, = \, c \odot x \, \oplus \, d \odot y $.
Let $i$ be an index which minimizes $\,a_i + z_i$.
We need to show that this minimum is attained at least twice.
By definition, $\,z_i \,$ is equal to either $c+x_i$ or $d+y_i$,
and, after permuting $x$ and $y$, we may assume
$\,z_i \, = \, c + x_i \, \leq \, d + y_i $.
 Since, for all $k$, $a_i+z_i \le a_k+z_k$ and $z_k \le c+x_{k}$, it follows
 that $a_i+x_{i}\le a_k+x_{k}$ for all $k$, so that $a_i+x_{i}$ achieves
the minimum of $\, \{a_1+x_1, \ldots,a_n+x_n\}$.
Since $x$ is in $H$, there exists some index $j \not= i$ for which
$a_i+x_i=a_j+x_j$.
But now $\,a_j + z_j \le a_j + c + x_j = c + a_i + x_{i} = a_i + z_i$.
 Since $a_i+z_i$ is the minimum of all
$a_j+z_j$, the two must be equal, and this minimum is obtained at least twice as desired.
\end{proof}

Proposition \ref{hypconv} implies that if $V$
is a subset of $\TP^{n-1}$ which happens to lie
in a tropical hyperplane $H$, then its
tropical convex hull $\,{\rm tconv}(V)\,$ will lie in $H$ as well.
The same holds for tropical planes of higher codimension.
Recall that every \emph{tropical plane}
 is an intersection of tropical hyperplanes \cite{SS}.
But the converse does not hold: not every intersection of tropical
hyperplanes qualifies as a tropical plane (see \cite[\S 5]{RGST}).
Proposition  \ref{hypconv}
and the first statement in Theorem
\ref{LotsOfStuff} imply:

\begin{cor} \label{planesconvex}
Tropical planes in $\TP^{n-1}$ are tropically convex.
\end{cor}

A theorem in classical geometry states
that every point outside a closed convex set
can be separated from the convex set by a hyperplane.
The same statement holds in tropical geometry.
This follows from the results in \cite{CGQ}.
Some caution is needed, however, since
the definition of hyperplane in \cite{CGQ}
differs from our definition of hyperplane,
as explained in \cite{RGST}. In our definition,
a tropical hyperplane is a fan which divides
$\TP^{n-1}$ into $n$ convex cones,
each of which is also tropically convex.
Rather than stating the most general
separation theorem, we will now
focus our attention on tropical polytopes,
in which case the separation theorem is
the Farkas Lemma stated in the next section.

\section{Tropical polytopes and cell complexes}\label{polytopes}

Throughout this section we fix a finite subset
$\, V = \{v_1,v_2,\ldots,v_r \} \,$ of tropical
projective space $\TP^{n-1}$. Here
$\,v_i = (v_{i1},v_{i2},\ldots,v_{in})$.
Our goal is to
study the tropical polytope $\, P =  {\rm tconv}(V)$.
We begin by describing the natural cell decomposition
of $\TP^{n-1}$ induced by the fixed finite subset $V$.

Let $x$ be any point in $\TP^{n-1}$.
The \emph{type} of $x$ relative to $V$ is the ordered $n$-tuple
$(S_1,\ldots,S_n)$ of subsets $S_j \subseteq \{1,2,\ldots,r\}$
which is defined as follows: An index $i$ is in $S_j$ if
$$ v_{ij} - x_j \quad = \quad {\rm min}
(  v_{i1} - x_1,  v_{i2} - x_2, \ldots,  v_{in} - x_n ). $$
Equivalently, if we set  $\,\lambda_i
= \, {\rm min}\{ \, \lambda \in \RR \, : \, \lambda \odot v_i \, \oplus x \, = \, x \,\}\,$
then $S_j$ is the set of all indices $i$ such that
$ \lambda_i \odot v_i \,$ and $x$ have the same $j$-th coordinate. We say that 
an $n$-tuple of indices $S=(S_1,\ldots,S_n)$ is a \emph{type} if it arises in this manner. Note that 
every $i$ must be in some $S_j$.

\begin{ex} \rm
Let $r=n=3$, $v_1 = (0,0,2),\,
v_2 = (0,2,0)$ and $v_3 = (0,1,-2)$.
There are $30$ possible types as $x$ ranges over the
plane $\TP^2$. The corresponding cell decomposition
has six convex regions (one bounded, five unbounded),
$15$ edges ($6$ bounded, $9$ unbounded)
and $6$ vertices. For instance, the point
$x = (0,1,-1)$ has $\, {\rm type}(x) = \bigl( \{2\}, \{1\}, \{3\} \bigr)\,$
and its cell is a bounded pentagon.
The point $ x'= (0,0,0)$ has
$\, {\rm type}(x') =  \bigl( \{1,2\}, \{1\}, \{2,3\} \bigr)\,$
and its cell is one of the six vertices.
The point $x'' = (0,0,-3)$ has
$\, {\rm type}(x'') = \bigl\{ \{1,2,3\}, \{1\}, \emptyset \bigr)\,$
and its cell is an unbounded edge.
\end{ex}

Our first application of types is the following separation theorem.

\begin{prop}[Tropical Farkas Lemma]\label{farkas}
For all $x \in \TP^{n-1}$, exactly one of the following is true.

(i) the point $x$ is in the tropical polytope $P = {\rm tconv}(V)$, or

(ii) there exists a tropical hyperplane which separates $x$ from $P$.
\end{prop}

The separation statement in
part (ii) means the following: if the hyperplane
is given by  (\ref{hyperplane}) and
$\, a_k + x_k \, = \, {\rm min}(a_1+x_1,\ldots,a_n+x_n)\,$
then  $\,a_k + y_k \, > \,  {\rm min}(a_1+y_1,\ldots,a_n+y_n)\,$
for all $y \in P$.

\begin{proof}
Consider any point $x \in \TP^{n-1}$, with $\,
{\rm type}(x)=(S_1,\ldots,S_n)$,	
and let $\,\lambda_i
= \, {\rm min}\{ \, \lambda \in \RR \,
 : \, \lambda \odot v_i \, \oplus x \, = \, x \,\}\,$
as before.
We define
\begin{equation}
\label{nearestpointmap}
\pi_V(x) \quad = \quad \lambda_1 \odot v_1 \,
\oplus \, \lambda_2 \odot v_2 \, \oplus \,\cdots \, \oplus \,
\lambda_r \odot v_r  .
\end{equation}
There are two cases: either $\,\pi_V(x) = x \,$ or $\,\pi_V(x) \not= x$.
The first case implies (i). Since (i) and (ii) clearly cannot occur
at the same time, it suffices to prove that the second case implies (ii).

Suppose that $\,\pi_V(x) \not= x$. Then
$S_k$ is empty for some index $k \in \{1,\ldots,n\}$. This means that
 $\,v_{ik}+\lambda_i-x_k>0 \,$ for $i=1,2,\ldots,r$. Let
 $\varepsilon > 0$ be smaller than any of these $r$ positive reals.
We now choose our separating tropical hyperplane   (\ref{hyperplane}) as follows:
\begin{equation}
\label{separator}
a_k \,\, := \,\, - x_k - \varepsilon \qquad \hbox{and} \qquad
a_j \,\,:= \,\,-x_j \, \, \,\,\hbox{for} \,\,\, j \in \{1,\ldots,n\}\backslash \{k\}.
\end{equation}
This certainly satisfies $a_k+x_k = {\rm min}(a_1+x_1.\ldots,a_n+x_n)$.
Now, consider any point $\,y=\bigoplus_{i=1}^r c_i\odot v_i \,$
in ${\rm tconv}(V)$.  Pick any $m$ such that $y_k = c_m + v_{mk}$. By
definition of the $\lambda_i$,  we have 
$\,x_k \leq \lambda_m + v_{mk}\,$ for all $k$, 
and there exists some $j$ with
$\,x_j = \lambda_m + v_{mj}$. These equations
and inequalities imply
\begin{eqnarray*}
& a_k + y_k \quad
= \,\, \,
a_k + c_m + v_{mk} 
\,\ =  \,\,c_m + v_{mk} - x_k - \varepsilon 
\,\ >  \,\,   
 c_m - \lambda_m \\ &
\,\, = \,\,c_m + v_{mj} - x_j
\,\, \geq \,\, y_j-x_j 
\,\, = \,\, a_j+y_j
\quad \geq  \quad
{\rm min}(a_1+y_1,\ldots,a_n+y_n).
\end{eqnarray*}
 Therefore, the hyperplane defined by (\ref{separator}) separates $x$ from $P$ as desired.
\end{proof}

The construction in (\ref{nearestpointmap})
defines a map $\,\pi_V \,: \, \TP^{n-1} \rightarrow P \,$
whose restriction to $P$ is the identity.
This map is the tropical version of the
\emph{nearest point map} onto a closed convex set
in ordinary geometry. Such maps were studied
in \cite{CGQ} for convex subsets in arbitrary idempotent semimodules.

If $S = (S_1,\ldots,S_n)$ 
and $T = (T_1,\ldots,T_n)$ are $n$-tuples
of subsets of $\{1,2,\ldots,r\}$, then we write
$\, S \subseteq T \,$ if $\, S_j \subseteq T_j \,$
for $j=1,\ldots,n$. We also
consider the set of all points whose type contains $S$:
$$ X_S \quad := \quad \bigl\{\,
x \in \TP^{n-1} \, : \, S \,\subseteq \,{\rm type}(x)\bigr\}. $$

\begin{lem}\label{ineqs}
The set $X_S$ is a closed convex polyhedron
(in the usual sense). More precisely,
\begin{equation}
\label{X_Spoly}
 X_S \quad = \quad \bigl\{\,
x \in \TP^{n-1} \,:\,
 x_k-x_j \le v_{ik}-v_{ij}  \,\,\,\,
\hbox{for all} \,\, 
j,k \in \{ 1,\ldots,n\} \,\,
\hbox{such that} \,\, i \in S_j \,\bigr\} .
\end{equation}
\end{lem}

\begin{proof}
Let $x \in \TP^{n-1}$ and $T = {\rm type}(x)$.
First, suppose $x$ is in $X_S$. Then $S \subseteq T$.
For every $i,j,k$ such that 
$i\in 
S_j$, we also have $i\in T_j$, and so by definition we have $v_{ij}-x_j \le v_{ik}-x_k$, or 
$x_k-x_j\le v_{ik}-v_{ij}$. Hence $x$ lies in the
set on the right hand side of (\ref{X_Spoly}).
For the proof of the reverse inclusion,
suppose that $x$ lies in the right hand side of 
(\ref{X_Spoly}).
Then, for all $i,j$ with $i\in S_j$, and for all $k$, we have  
$v_{ij}-x_j\le v_{ik}-x_k$. This means that
$\,v_{ij}-x_j =\,{\rm  min}(v_{i1}-x_1,\ldots,v_{in}-x_n)\,$
and hence $\,i\in T_j$. Consequently, for 
all $j$, we have $S_j\subset T_j$, and so $x\in X_S$.
\end{proof}

\begin{figure}
 \input{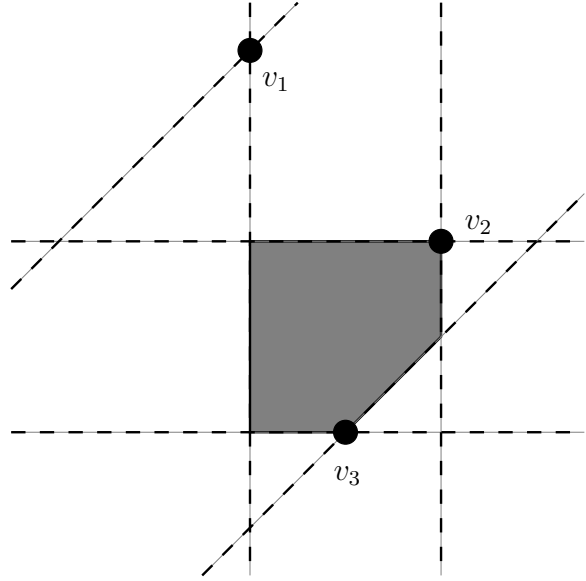}
\caption{\label{reg213}
The region $X_{(2,1,3)}$ in the tropical convex hull of 
$v_1$, $v_2$ and $v_3$.}
\end{figure}

As an example for Lemma~\ref{ineqs}, we 
consider the region $X_{(2,1,3)}$ in
the tropical convex hull of $v_1=(0,0,2)$, $v_2=(0,2,0)$, and $v_3=(0,1,-2)$. 
This region is defined by six linear inequalities, one of which is redundant, as depicted in Figure~\ref{reg213}.
Lemma~\ref{ineqs} has the following 
immediate corollaries.

\begin{cor}\label{intersectregion}
The intersection $X_S\cap X_T$ is equal to the polyhedron $\,X_{S\cup T}$.
\end{cor}
\begin{proof}
The inequalities defining $X_{S\cup T}$ are precisely the union of the inequalities defining $X_S$ 
and $X_T$, and points satisfying these inequalities are precisely those in $X_S\cap X_T$.
\end{proof}

\begin{cor}\label{boundedness}
The polyhedron $X_S$ is bounded if and only if
 $S_j \not= \emptyset $ for all $j = 1,2,\ldots,n$.
\end{cor}

\begin{proof}
Suppose that $S_j\neq \emptyset$ for all $j=1,2,\ldots,n$. Then for every $j$ and $k$, we can find 
$i\in S_j$ and $m\in S_k$, which via Lemma~\ref{ineqs} yield the inequalities 
$v_{mk}-v_{mj}\le x_k-x_j\le 
v_{ik}-v_{ij}$. This implies that each $x_k-x_j$ is bounded on $X_S$, which means that $X_S$ is 
a bounded subset of $\TP^{n-1}$.

Conversely, suppose some $S_j$ is empty. Then the only inequalities involving $x_j$ are of the 
form $x_j-x_k\le c_{jk}$. Consequently, if any point $x$ is in $S_j$, so too is $x-ke_j$ for 
$k>0$, where $e_j$ is the $j$-th basis vector. Therefore, in this case, $X_S$ is unbounded.
\end{proof}

\begin{cor}\label{regionfaces}
Suppose we have $S=(S_1,\ldots,S_n)$, with $S_1\cup\cdots\cup S_n = \{1,\ldots,r\}$. 
Then if $S\subseteq T$, $X_T$ is a
face of $X_S$, and furthermore all faces of $X_S$ are of this form. 
\end{cor}

\begin{proof}
For the first part, it suffices to prove that the statement is true when $T$ covers $S$ in the
poset of containment, i.e. when $T_j=S_j\cup \{i\}$ for 
some $j\in \{1,\ldots,n\}$ and $i\not\in
S_j$, and $T_k=S_k$ for $k\neq j$.

We have the inequality presentation of $X_S$ given by Lemma~\ref{ineqs}. By the same lemma, the
inequality presentation of $X_T$ consists of the 
inequalities defining $X_S$ together
with the inequalities 
\begin{equation}\label{extraineqs}
\{x_k-x_j \le v_{ik}-v_{ij}\,\mid\, k\in \{1,\ldots,n\}\}.
\end{equation}

By assumption, $i$ is in some $S_m$. We claim that $X_T$ is the face of $S$ defined by the equality
\begin{equation}\label{extraeq} x_m-x_j = v_{im}-v_{ij}. \end{equation}
Since $X_S$ satisfies the inequality $x_j-x_m\le v_{ij}-v_{im}$, (\ref{extraeq}) defines a face $F$ of 
$S$. The inequality $x_m-x_j\le v_{im}-v_{ij}$ is in the set (\ref{extraineqs}), so (\ref{extraeq}) is valid on 
$X_T$ and $X_T\subseteq F$. However, any point in $F$, being in $X_S$, satisfies $x_k-x_m\le v_{ik}-v_{im}$ for all 
$k\in \{1,\ldots,n\}$. Adding (\ref{extraeq}) to these inequalities proves that the inequalities 
(\ref{extraineqs}) are valid on $F$, and hence $F\subseteq X_T$. So $X_T=F$ as desired.

By the discussion
in the proof of the first part, prescribing equality in the facet-defining inequality $x_k-x_j\le v_{ik}-v_{ij}$ 
yields $X_T$, where $T_k = S_k\cup \{i\}$ and $T_j=S_j$ for $j\neq k$. Therefore, all facets of $X_S$ can be obtained 
as regions $X_T$, and it follows recursively that all faces of $X_S$ are of this form.
\end{proof}

\begin{cor}\label{alwaystype}
Suppose that $S=(S_1,\ldots,S_n)$ is an $n$-tuple of indices satisfying 
$\,S_1\cup\cdots\cup S_n=\{1,\ldots,r\}$.
Then $X_S$ is equal to $X_T$ for some type $T$.   
\end{cor}

\begin{proof}
Let $x$ be a point in the relative interior of $X_S$, 
and let $\,T = {\rm type}(x)$. Since $x\in
X_S$, $T$ contains $S$, and by Lemma~\ref{regionfaces}, 
$X_T$ is a face of $X_S$.
However, since $x$ is in the relative interior of $X_S$, 
the only face of $X_S$ containing $x$ is
$X_S$ itself, so we must have $X_S=X_T$ as desired.
\end{proof}

We are now prepared to state our main theorem in this section.

\begin{thm}\label{celldecomp}
The collection of convex polyhedra $X_S$, where $S$ ranges over all types, defines
a cell decomposition $ \,{\mathcal C}_V$ of $\,\TP^{n-1}$.
The tropical polytope $\, P = {\rm tconv}(V)\,$ equals
the union of all bounded cells $X_S$ in this decomposition.
\end{thm}

\begin{proof}
Since each point has a type, it is clear that the union of the $X_S$ is 
equal to $\TP^{n-1}$. By Corollary~\ref{regionfaces}, the faces of $X_S$ are equal to $X_U$ for $S\subseteq U$, and 
by Corollary~\ref{alwaystype}, $X_U=X_W$ for some type $W$, and hence $X_U$ is in our collection. The only thing 
remaining to check to show that this collection defines a cell decomposition is that $X_S\cap X_T$ is a face of both 
$X_S$ and $X_T$, but $X_S\cap X_T=X_{S\cup T}$ by Corollary~\ref{intersectregion}, and $X_{S\cup T}$ is a face of 
$X_S$ and $X_T$ by Corollary~\ref{regionfaces}.

For the second assertion consider any point
$x \in \TP^{n-1}$ and let $S = {\rm type}(x)$.
We have seen in the proof of the
Tropical Farkas Lemma (Proposition~\ref{farkas})
that $x $ lies in $P$ if and only if no $S_j $ is empty.
By Corollary~\ref{boundedness}, this is equivalent to
the polyhedron $X_S$ being bounded.
\end{proof}

The collection of bounded cells $X_S$ is referred to as the tropical complex generated by $V$; thus, 
Theorem~\ref{celldecomp} states that this provides a polyhedral decomposition of the polytope $P={\rm 
tconv}(V)$. 
Different sets $V$ may 
have the same tropical polytope as their convex hull, but generate different tropical complexes; the decomposition 
of a tropical polytope depends on the chosen generating set, although we will see later 
(Proposition~\ref{vertexunique}) that there is a unique minimal generating set.

Here is a nice geometric construction of
the cell decomposition $\, {\mathcal C}_V \,$
of $\TP^{n-1}$ induced by $\, V = \{v_1,\ldots,v_r\}$.
Let ${\mathcal F}$ be the fan in $\TP^{n-1}$
defined by the tropical hyperplane (\ref{hyperplane})
with  $\,a_1 = \cdots = a_n = 0$. 
Two vectors $x$ and $y$ lie in the same
relatively open cone of the fan ${\mathcal F}$
if and only if
$$ \{\, j \,: \, x_j \, = \, {\rm min}(x_1,\ldots,x_n) \,\}
\quad = \quad 
\{\, j \,: \, x_j \, = \, {\rm min}(y_1,\ldots,y_n) \,\}. $$
If we translate the negative of ${\mathcal F}$
by the vector $v_i$ then we get a new fan which we denote by
$\, v_i - {\mathcal F} $. 
Two vectors $x$ and $y$ lie in the same
relatively open cone of the fan  $\, v_i - {\mathcal F} \,$
if and only
$$ \{\, j \,: \, x_j-v_{ij}\, = \, {\rm max}(x_1-v_{i1},\ldots,x_n-v_{in}) \,\}
\quad = \quad 
\{\, j \,: \, y_j-v_{ij} \, = \, {\rm max}(y_1-v_{i1},\ldots,y_n-v_{in}) \,\}. $$

\begin{prop}\label{fans}
The cell decomposition $\, {\mathcal C}_V\,$ is the common refinement of
the $r$ fans $\, v_i - {\mathcal F}$.
\end{prop}

\begin{proof}
We need to show that the cells of this common refinement are 
precisely the convex polyhedra $X_S$. Take 
a point $x$, with $\, T = {\rm type}(x)\,$
and define $S_x=(S_{x1},\ldots,S_{xn})$ by letting $i\in S_{xj}$ whenever
\begin{equation}
\label{maxnotmin}
x_j -v_{ij}\quad = \quad {\rm max}(x_1-v_{i1},\ldots,x_n-v_{in}).
\end{equation}
Two points $x$ and $y$ are in the relative interior of the same 
cell of the common refinement if 
and only if they are in 
the same relatively open cone of each fan; this is tantamount to saying that $S_x=S_y$. However, 
we claim that 
$S_x = T$. Indeed, taking the negative of 
both sides of (\ref{maxnotmin})
yields exactly the condition for $i$ being in $T_j$,
 by the definition of type. 
Consequently, the condition for two points having the same type is the same as the condition for 
the two points being in the relative interior of the same 
chamber of the common refinement of the fans
$\,  v_1 - {\mathcal F},\,
v_2 - {\mathcal F},\,\ldots,v_r - {\mathcal F}$.
\end{proof}

\begin{figure}
\input{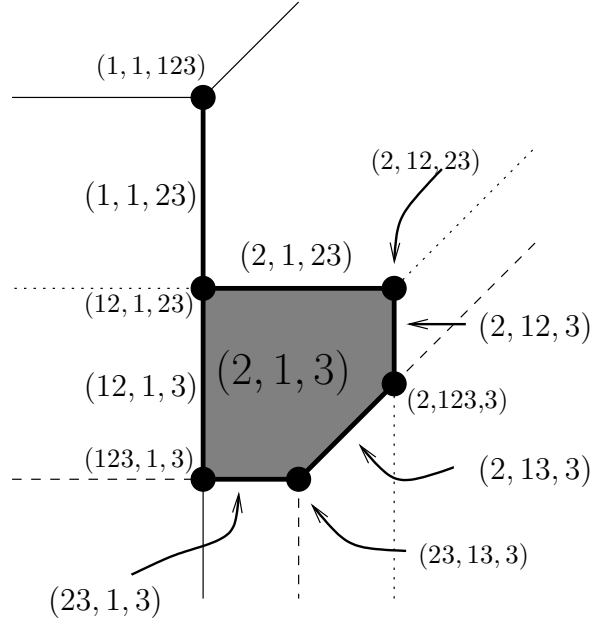}
\caption{\label{tropfans}
A tropical complex 
expressed as the bounded cells in the common refinement
of the fans  
$ v_1 - {\mathcal F}$,
$ v_2 - {\mathcal F}$ and
$ v_3 - {\mathcal F}$.
Cells are labeled with their types.}
\end{figure}

An example of this construction is shown for our usual example, where $v_1=(0,0,2)$, $v_2=(0,2,0)$, and 
$v_3=(0,1,-2)$, in Figure~\ref{tropfans}.

The next few results provide additional information
about the polyhedron $X_S$. 
Let $G_S$ denote the  undirected  
graph with vertices $\{1,\ldots,n\}$, 
where  $\{j,k\}$ is an edge
if and only if $\,S_j \cap S_k\neq \emptyset$.

\begin{prop}\label{dimension}
The dimension $d$ of the polyhedron $X_S$ is one less than the number of 
connected components of $G_S$, and $X_S$ 
is affinely and tropically isomorphic to some polyhedron
$X_T$ in $\TP^{d}$.
\end{prop} 

\begin{proof}
The proof is by induction on $n$. 
Suppose we have $i\in S_j\cap S_k$. Then $X_S$ satisfies the linear
equation $x_k-x_j = c$ where $c=v_{ik}-v_{ij}$. 
Eliminating the variable $x_k$ (projecting
onto $\TP^{n-2}$), $X_S$ is affinely and tropically isomorphic to 
$X_T$ where the type $T$ is defined by  $T_r  = S_r$ for
$r\neq j$ and $T_j = S_j\cup S_k$. The region $X_T$ exists in the 
cell decomposition of $\TP^{n-2}$ induced by the vectors
$w_1,\ldots,w_n$ with  $w_{ir} = v_{ir}$ for $r\neq j$, 
and $w_{ij} = \text{max}(v_{ij}, v_{ik}-c)$. 
The graph $G_T$ is 
obtained from the graph $G_S$ by contracting the edge 
$\{j,k\}$, and thus has the same number of 
connected components.

This induction on $n$ reduces us to the case where all of the $S_j$ are pairwise disjoint. We must show that  $X_S$ has dimension $n-1$.
Suppose not. Then $X_S$ lies in $\TP^{n-1}$ but has dimension less than
$n-1$. Therefore, one of the inequalities in  (\ref{X_Spoly}) holds with 
equality, say
$\, x_k - x_j \, = \, v_{ik} - v_{ij} \,$ for all $x \in X_S$. The inequality
``$\leq $'' implies $i \in S_j$ and the inequality
``$\geq $'' implies $i \in S_k$. Hence $S_j$ and $S_k$
are not disjoint, a contradiction.
\end{proof}

The following proposition can be regarded as a converse to
Lemma~\ref{ineqs}.

\begin{prop}\label{convineqs}
Let $R$ be any polytope in $\TP^{n-1}$ 
defined by inequalities of the form $\,x_k-x_j\le c_{jk}$.
Then $R$ arises as a cell $X_S$ in the decomposition ${\mathcal C}_V$
of  $\TP^{n-1}$ defined by some set $V = \{v_1,\ldots,v_n\}$.
\end{prop}

\begin{proof}
Define the vectors $v_i$ to have coordinates $v_{ij} = c_{ij}$ for $i\neq j$, 
and $v_{ii} = 0$. (If $c_{ij}$ did not appear in the given inequality
presentation then simply take it to be a very large positive number.)
Then by Lemma~\ref{ineqs}, the polytope in $\TP^{n-1}$ defined by
the inequalities $\, x_k-x_j\le c_{jk} \,$ is
precisely the unique cell of type $(1,2,\ldots,n)$ in the tropical 
conver hull of  $\{v_1,\ldots,v_n\}$.
\end{proof}

The region $X_S$ is a polytope both in the
ordinary sense and in the tropical sense.

\begin{prop}\label{regconvex} 
Every
bounded cell $X_S$ in the tropical complex generated by $V$ is itself a tropical polytope,
equal to the tropical convex hull of its vertices. 
The number of vertices of the polytope $X_S$ is at most $\, \binom{2n-2}{n-1}\,  $,
and this bound is tight for all positive integers $n$.
\end{prop}

\begin{proof}
By Proposition~\ref{dimension}, if $X_S$ has dimension $d$, it is affinely and tropically
isomorphic to a region in the convex hull of a set of points in $\TP^{d}$, so it suffices to 
consider the full-dimensional case.

The inequality presentation of Lemma~\ref{ineqs} demonstrates that $X_S$ is tropically convex for 
all $S$, since if two points satisfy an inequality of that form, so does 
any tropical linear combination 
thereof. Therefore, it suffices to show that $X_S$ is contained in the tropical convex hull of its 
vertices.

The proof is by induction on the dimension of $X_S$. All proper faces of $X_S$ are 
polytopes
$X_T$ of lower dimension, and by induction are contained in the tropical convex hull of their 
vertices. These vertices are a subset of the vertices of $X_S$, and so this face is
in the tropical convex hull.

Take any point $x=(x_1,\ldots,x_n)$ in the interior of $X_S$. Since $X_S$ has dimension
$n$, we can travel in any direction from $x$ while remaining in $X_S$. Let us travel in the
$(1,0,\ldots,0)$ direction until we hit the boundary, to obtain points
$y_1=(x_1+b,x_2,\ldots,x_n)$ and $y_2=(x_1-c,x_2,\ldots,x_n)$ in the boundary of $X_S$. These
points are contained in the tropical convex hull by the inductive hypothesis, which means that $x
= y_1\oplus c\odot y_2$ is also, completing the proof of the first assertion.

For the second assertion, we consider the convex hull
of all differences of unit vectors, $e_i - e_j$. This is a lattice polytope
of dimension $n-1$ and normalized volume $\,\binom{2n-2}{n-1}$. To see this, we observe that this polytope is tiled
by $n$ copies of the convex hull of the origin and 
the $\,\binom{n}{2}$ vectors $e_i-e_j$ with $i<j$. 
The other $n-1$ copies are gotten by
 cyclic permutation of the coordinates.
This latter polytope was studied by
Gel'fand, Graev and Postnikov, who showed in \cite[Theorem 2.3 (2)]{GGP} that the normalized volume of this polytope 
equals the Catalan number $\frac{1}{n} \binom{2n-2}{n-1}$.

We conclude that every complete fan whose rays are among the
vectors $e_i - e_j$ has at most $\, \binom{2n-2}{n-1}  \,$ maximal cones.
This applies in particular to the normal fan of $X_S$, 
hence $X_S$ has at most $\, \binom{2n-2}{n-1}  \,$ vertices.
Since the configuration $\{e_i-e_j\}$ is unimodular,
the bound is tight whenever the fan is simplicial and uses all the
rays $e_i-e_j$.
\end{proof}

We close this section with two more results
about arbitrary tropical polytopes in $\TP^{n-1}$.

\begin{prop}\label{intersectpolyt}
If $P$ and $Q$ are tropical polytopes in $\TP^{n-1}$ then
 $P\cap Q$ is also a tropical polytope. 
\end{prop}

\begin{proof}
Since $P$ and $Q$ are both tropically convex, $P\cap Q$ must also be. Consequently, if we can find 
a finite set of points in $P\cap Q$ whose convex hull contains all of $P\cap Q$, we will be done.  By Theorem \ref{celldecomp},
$P$ and $Q$ are the finite unions of
bounded cells  $\{X_S\}$ and $\{X_T\}$ respectively, so $P\cap 
Q$ is the finite union of the cells
$X_S\cap X_T$. 
Consider any $X_S\cap X_T$. Using Lemma~\ref{ineqs} to obtain the inequality representations of 
$X_S$ and $X_T$, we see that this region is of the form dictated by Proposition~\ref{convineqs}, 
and therefore obtainable as a cell $X_W$ in some tropical complex. By Proposition~\ref{regconvex}, 
$X_W$ is itself a tropical polytope, and we can therefore find a finite set whose convex hull is 
equal to $X_S\cap X_T$. Taking the union of these sets over all choices of $S$ and $T$ then gives 
us the desired set of points whose convex hull contains all of $P\cap Q$.
\end{proof}

\begin{prop}\label{vertexunique}
Let $P\subset \TP^{n-1}$ be a tropical polytope. Then there exists a unique minimal set $V$ such that $P={\rm 
tconv}(V)$.
\end{prop}

\begin{proof}
Suppose that $P$ has two minimal generating sets, 
$V=\{v_1,\ldots,v_m\}$ and $W=\{w_1,\ldots,w_r\}$. 
Write each element of $W$ as $w_i = \oplus_{j=1}^m c_{ij}\odot v_j$.
We claim that $V \subseteq W$.
Consider $v_1\in V$ and write
\begin{equation}\label{selfexpr}
v_1 \,\,\,= \,\,\, \bigoplus_{i=1}^r d_i \odot w_i 
\,\, \,\, = \,\,\,\bigoplus_{j=1}^m f_j\odot v_j
\qquad
\hbox{where $f_j = {\rm min}_i (d_i + c_{ij})$}.
\end{equation}
 If the term $f_1\odot v_1$ does not minimize any coordinate in the 
right-hand  side of (\ref{selfexpr}), 
then $v_1$ is a linear combination of $v_2,\ldots,v_m$, 
contradicting the minimality of $V$. 
However, if $f_1\odot v_1$ minimizes any coordinate in this expression, 
it must minimize all of them, since 
$(v_1)_j-(v_1)_k = (f_1\odot v_1)_j - (f_1\odot v_1)_k$. 
In this case we get $\,v_1 
= f_1 \odot v_1$, or $f_1 = 0$. Pick any $i$ 
for which $f_1 = d_i + c_{i1}$; we claim that $w_i = c_{i1} \odot v_1$. 
Indeed, if any other term in $w_i = 
\oplus_{j=1}^m c_{ij}\odot v_j$ contributed nontrivially to $w_i$, 
that term would also contribute  to 
the expression on the right-hand side of (\ref{selfexpr}), 
which is a contradiction.
Consequently, $V\subseteq W$, which means $V=W$ since 
both sets are minimal by hypothesis.
\end{proof}

Like many of the results presented in this section, 
Propositions~\ref{intersectpolyt} and \ref{vertexunique} 
parallel results on ordinary polytopes. 
We have already mentioned the tropical analogues of the Farkas 
Lemma and of Carath\'{e}odory's 
Theorem (Propositions~\ref{carat} and \ref{farkas}); 
Proposition~\ref{dimension} is analogous to the result that a 
polytope $P\subset \RR^n$ of dimension $d$ is 
affinely isomorphic to some $Q\subset \RR^d$. 
Proposition~\ref{regconvex} hints at a duality between an inequality 
representation and a vertex representation of a tropical polytope; this duality has been studied in greater detail 
by Michael Joswig~\cite{Jos}.

\section{Subdividing products of simplices}\label{triang}

Every set $V=\{v_1,\ldots,v_r\}$ of $r$ points in $\TP^{n-1}$ begets a tropical polytope
$P=\text{tconv}(V)$ equipped with a cell decomposition into the tropical complex generated by $V$.  Each
cell of this tropical complex is labelled by its type, which is an $n$-vector of finite subsets of
$\{1,\ldots,r\}$. Two configurations (and their corresponding tropical complexes) $V$ and $W$ have the
same \emph{combinatorial type} if the types occurring in their tropical complexes are identical; note 
that by
Lemma~\ref{regionfaces}, this implies that the face posets of these polyhedral complexes are isomorphic.

With the definition in the previous paragraph, the statement of Theorem \ref{main} has now finally been
made precise.  We will prove this correspondence between tropical complexes and subdivisions of products
of simplices by constructing the polyhedral complex ${\mathcal C}_P$ in a higher-dimensional space.

Let $W$ denote the $(r+n-1)$-dimensional real vector space 
$\, \RR^{r+n}/(1,\ldots,1,-1,\ldots,-1)$.
The natural coordinates on $W$ are denoted
$(y,z) = (y_1,\ldots,y_r, z_1,\ldots,z_n)$.
As before, we fix an ordered subset 
$V = \{v_1,\ldots,v_r\}$ of $\TP^{n-1}$
where $v_i = (v_{i1},\ldots,v_{in})$.
This defines the unbounded polyhedron
\begin{equation}
\label{masterpolyhedron}
{\mathcal P}_V \quad = \quad
\bigl\{\, (y,z) \in W \,:\,\, y_i + z_j \leq v_{ij} \,\,\,
\hbox{for all} \,\,i \in \{1,\ldots,r\}
\,\,\hbox{and}\,\, j \in \{1,\ldots,n\}\, \bigr\}.
\end{equation}

\begin{lem}\label{triangtrop}
There is a piecewise-linear isomorphism between
the tropical complex generated by $V$ and
the complex of bounded faces of the
$(r+n-1)$-dimensional polyhedron
$\,{\mathcal P}_V$. The image of a
cell $\,X_S\,$ of ${\mathcal C}_P$  under this isomorphism
is the
bounded face $\{y_i+z_j=v_{ij}\,:\,i\in S_j\}$ of
the polyhedron $\,{\mathcal P}_V$.
That bounded face  maps isomorphically to $\,X_S\,$ via 
projection onto the 
$z$-coordinates.
\end{lem}

\begin{proof}

Let $F$ be a bounded face of ${\mathcal P}_V$, and define $S_j$ via $i\in S_j$ if $y_i+z_j=v_{ij}$ is
valid on all of $F$. If some $y_i$ or $z_j$ appears in no equality, then we can subtract arbitrary
positive multiples of that basis vector to obtain elements of $F$, contradicting the assumption that $F$
is bounded. Therefore, each $i$ must appear in some $S_j$, and each $S_j$ must be nonempty.

Since every $y_i$ appears in some equality, given 
a specific $z$ in the projection of $F$ onto the
$z$-coordinates, there exists a unique $y$ for which $(y,z)\in F$, so this projection is an affine isomorphism from
$F$ to its image. We need to show that this image is equal to $X_S$. 

Let $z$ be a point in the image of this projection, coming from a
point $(y,z)$ in the relative interior of $F$. We claim that $z\in X_S$. 
Indeed, looking at the $j$th coordinate of $z$, we find
\begin{equation}\label{outtype}
 -y_i + v_{ij} \,\ge\, z_j \quad \hbox{for all} \,\, i, 
\end{equation}
\begin{equation}\label{intype}
-y_i + v_{ij} \,=\, z_j \quad \hbox{for} \,\,\,i\in S_j.
\end{equation}
The defining inequalities of $X_S$
are $x_j-x_k \le v_{ij}-v_{ik}$ with $i\in S_j$. Subtracting the
inequality $-y_i + v_{ik} \ge z_k$ from the equality in (\ref{intype}) yields that this inequality is valid on $z$ as
well.  Therefore, $z\in X_S$. Similar reasoning shows
that $\, S = {\rm type}(z)$.
We note that the relations (\ref{outtype}) and (\ref{intype}) 
can be rewritten elegantly in terms of
the tropical product of a row vector and a matrix:
\begin{equation}
\label{troplinalg}
z \quad = \quad (-y) \odot V \,\,\,= \,\,\,
\bigoplus_{i=1}^r (-y_i) \odot v_i.
\end{equation}

For the reverse inclusion, suppose that $\,z\in X_S$.
We define
 $\,y = V \odot (-z)$. This means that
\begin{equation}\label{expry}
y_i \,\,= \,\, {\rm min}
( v_{i1}-z_1, v_{i2}-z_2, \ldots, v_{in}-z_n).
\end{equation}
We claim that $(y,z)\in F$. Indeed, we certainly have $y_i + z_j \le v_{ij}$ for all $i$ and $j$, so $(y,z)\in
{\mathcal P}_V$. Furthermore, when $i\in S_j$, we know that $v_{ij}-z_j$ achieves the minimum in the right-hand side of
(\ref{expry}), so that $v_{ij}-z_j = y_i$ and $y_i+z_j = v_{ij}$ is satisfied. Consequently, $(y,z)\in F$ as desired.

It follows immediately that the two complexes are isomorphic: if $F$ is a face corresponding to $X_S$ and $G$ is a
face corresponding to $X_T$, where $S$ and $T$ are both types, then $X_S$ is a face of $X_T$ if and only if
$T\subseteq S$. However, by the discussion above, this is equivalent to saying that the equalities $G$ satisfies
(which correspond to $T$) are a subset of the equalities $F$ satisfies (which correspond to $S$); this is true if and 
only if $F$ is a face of $G$. So $X_S$ is a face of $X_T$ if and only if $F$ is a face of $G$, which implies the 
isomorphism of complexes.
\end{proof}

The boundary complex of the polyhedron $\,{\mathcal P}_V \,$  is polar 
to the regular subdivision of the 
product of simplices $\Delta_{r-1}\times \Delta_{n-1}$ 
defined by the weights $v_{ij}$. We denote this
regular polyhedral subdivision  
by  $\,(\partial {\mathcal P}_V)^*$. Explicitly,
a subset of vertices $(e_i,e_j)$
of $\Delta_{r-1}\times \Delta_{n-1}$ forms a cell
of  $\,(\partial {\mathcal P}_V)^*\,$
if and only if the equations
$\, y_i + z_j = v_{ij} \,$ indexed by these vertices
specify a face of the polyhedron ${\mathcal P}_V$.
We refer to the book of
De Loera, Rambau and Santos~\cite{DRS} for basics
on polyhedral subdivisions. 

We now present the proof of the result stated in the introduction.
\vskip .1cm

\noindent {\sl Proof of Theorem \ref{main}: }
The poset of bounded faces of ${\mathcal P}_V$
is antiisomorphic to the poset of interior cells
of the subdivision $\,(\partial {\mathcal P}_V)^*\,$ of 
$\Delta_{r-1} \times \Delta_{n-1}$. Since every  
full-dimensional cell of $\,(\partial {\mathcal P}_V)^*\,$ is interior,
the subdivision is uniquely determined by its interior cells.
In other words, the combinatorial type of ${\mathcal P}_V$
is uniquely determined by the lists of facets containing each bounded face of
${\mathcal P}_V$. These lists are
precisely the types of regions in ${\mathcal C}_P$
by Lemma \ref{triangtrop}. This completes the proof.
\qed
\vskip .3cm

Theorem~\ref{main}, which establishes a bijection between the tropical complexes generated by $r$ points
in $\TP^{n-1}$ and the regular subdivisions of a product of simplices $\Delta_{r-1} \times \Delta_{n-1}$,
has many striking consequences. First of all, we can pick off the types present in a tropical complex 
simply by looking at the cells present in the corresponding regular subdivision. In particular, if we have 
an interior cell $M$, the corresponding type appearing in the tropical complex is defined via $S_j = 
\{i\in [n]\,\mid\,(j,i)\in M\}$. 

It is worth noting that via the Cayley Trick~\cite{San}, 
Theorem~\ref{main} is equivalent to saying that
tropical complexes generated by $r$ points in $\TP^{n-1}$ are in bijection with the regular mixed
subdivisions of the dilated simplex $r\Delta^{n-1}$. This connection is expanded upon and employed in a
paper with Francisco Santos~\cite{DSS}.
Another astonishing consequence of Theorem~\ref{main}
is the identification of
the row span and column span of a matrix. This result can also be derived
from \cite[Theorem 42]{CGQ}.

\begin{thm}\label{polar}

Given any matrix $M\in \RR^{r\times n}$, the tropical complex generated by 
its
column vectors is isomorphic to the tropical complex generated by its row
vectors. This isomorphism is gotten by restricting the piecewise linear
maps $ \, \RR^n \rightarrow \RR^r,\, z \mapsto M\odot (-z) \,$ and $ \,\,
\RR^r \rightarrow \RR^n,\, y \mapsto (-y) \odot M $.

\end{thm}

\begin{proof}
By Theorem~\ref{main}, the matrix $M$ corresponds via the polyhedron
${\mathcal P}_M$ to a regular subdivision of $\Delta_{r-1} \times
\Delta_{n-1}$, and the complex of interior faces of this regular
subdivision is combinatorially isomorphic to both the tropical complex 
generated by its row vectors, which are $r$ points in $\TP^{n-1}$, and 
the tropical complex generated by its 
column vectors, which are $n$ points in $\TP^{r-1}$. Furthermore,
Lemma~\ref{triangtrop} tells us that the cell in ${\mathcal P}_M$ is
affinely isomorphic to its corresponding cell in both tropical complexes.
Finally, in the proof of Lemma~\ref{triangtrop}, we showed that the
point $(y,z)$ in a bounded face $F$ of ${\mathcal P}_M$ satisfies $y =
M\odot (-z)$ and $z = (-y) \odot M$. This point projects to $y$ and $z$,
and so the piecewise-linear isomorphism mapping these two complexes to
each other is defined by the stated maps.  
\end{proof}

The common tropical complex of these two tropical polytopes is given by the complex of bounded faces of
the common polyhedron ${\mathcal P}_M$, which lives in a space of dimension $r+n-1$; the tropical
polytopes are unfoldings of this complex into dimensions $r-1$ and $n-1$. Theorem~\ref{polar} also gives a
natural bijection between the combinatorial types of tropical convex hulls of $r$ points in $\TP^{n-1}$
and the combinatorial types of tropical convex hulls of $n$ points in $\TP^{r-1}$, incidentally proving
that there are the same number of each. This duality statement extends a similar statement in~\cite{CGQ}.

\begin{figure}
\input{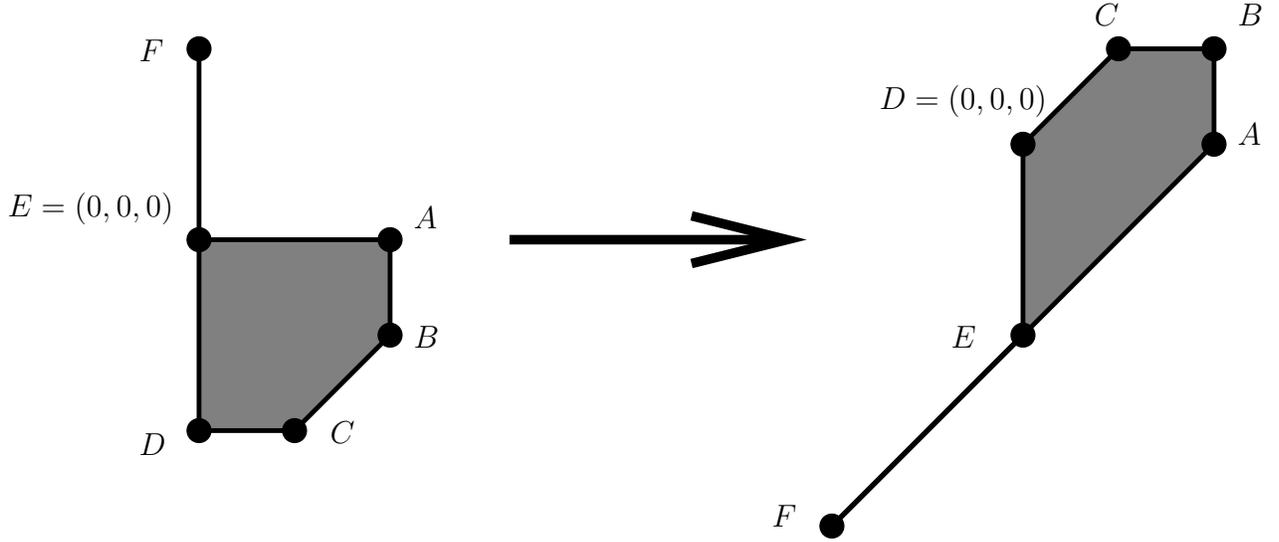}
\caption{\label{dualtriang}
A demonstration of tropical polytope duality.}
\end{figure}

Figure~\ref{dualtriang} shows the dual of the convex hull of $\{(0,0,2),(0,2,0),(0,1,-2)\}$, also a 
tropical triangle (here $r=n=3$). For instance, we compute:
\[
\left(
\begin{array}{ccc}
0 & 0 & 2 \\
0 & 2 & 0 \\
0 & 1 & -2 \\
\end{array}
\right)
\left(
\begin{array}{c}
0 \\
0 \\
-2 \\
\end{array}
\right)
= 
\left(
\begin{array}{c}
0 \\
-2 \\
-4 \\
\end{array}
\right).
\]
This point
is the image of the point $(0,0,2)$ under this duality map. Note that duality does not preserve the
generating set; the polytope on the right is the convex hull of points $\{F,D,B\}$, while the polytope on
the left is the convex hull of points $\{F,A,C\}$. This is necessary, of course, since in general a
polytope with $r$ vertices is mapped to a polytope with $n$ vertices, and $r$ need not equal $n$ as it
does in our example.

We now discuss the generic case when the subdivision
$\,(\partial {\mathcal P}_V)^*\,$ is a regular triangulation
of $\Delta_{r-1} \times \Delta_{n-1}$.
We refer to \cite[\S 5]{RGST}
for the geometric interpretation of the
\emph{tropical determinant}.

\begin{prop} \label{thegenericcase}
 For a configuration  $V$ of $r $ points in
$\TP^{n-1}$ with $r \geq n$ the following are equivalent:
\begin{enumerate}
\item The regular subdivision $(\partial {\mathcal P}_V)^*$ is 
 a triangulation of $\Delta_{r-1} \times \Delta_{n-1}$.
\item No $k$ of the points in $V$ have projections onto a $k$-dimensional coordinate subspace which lie in a tropical 
hyperplane, for any $2\le k\le n$. 
\item No $k \times k$-submatrix of the
$r \times n$-matrix $(v_{ij})$ is tropically singular, i.e. has vanishing tropical determinant, for any $2\le k\le 
n$.
\end{enumerate}
\end{prop}

\begin{proof}
The last equivalence is proven in~\cite[Lemma 5.1]{RGST}. We will prove that (1) and (3) are equivalent. The tropical determinant of
a  $k$ by $k$ matrix $M$ is the tropical polynomial
$\oplus_{\sigma\in S_k} (\odot_{i=1}^k M_{i\sigma(i)})$. 
The matrix $M$ is tropically singular
if the minimum ${\rm
min}_{\sigma\in S_k}\,(\sum_{i=1}^k M_{i\sigma(i)})$ is achieved twice.

The regular subdivision $(\partial {\mathcal P}_V)^*$ is a triangulation if and only if the polyhedron ${\mathcal P}_V$
is simple, which is to say if and only if no $r+n$ of the facets $y_i+z_j\le v_{ij}$ meet at a single vertex. For
each vertex $v$, consider the bipartite graph $G_v$ consisting of vertices $y_1,\ldots,y_n$ and $z_1,\ldots,z_j$ with
an edge connecting $y_i$ and $z_j$ if $v$ lies on the corresponding facet. This graph is connected, since each $y_i$
and $z_j$ appears in some such inequality, and thus it will have a cycle if and only if it has at least $r+n$ edges.
Consequently, ${\mathcal P}_V$ is not simple if and only there exists some $G_v$ with a cycle.

If there is a cycle, without loss of generality it reads $y_1, z_1, y_2,
z_2, \ldots, y_k, z_k$. Consider the submatrix $M$ of $(v_{ij})$ given by
$1\le i\le k$ and $1\le j\le k$. We have $y_1+z_1=M_{11}$,
$y_2+z_2=M_{22}$, and so on, and also $z_1+y_2=M_{12}, \ldots,z_k+y_1 =
M_{k1}$. Adding up all of these equalities yields
$y_1+\cdots+y_k+z_1+\cdots+z_k= M_{11}+\cdots+M_{kk} =
M_{12}+\cdots+M_{k1}$. But consider any permutation $\sigma$ in the 
symmetric group $S_k$. Since
we have $M_{i\sigma(i)}=v_{i\sigma(i)} \ge y_i+z_{\sigma(i)}$, we have
$\sum M_{i\sigma(i)}\ge x_1+\cdots+x_k+y_1+\cdots+y_k$. Consequently, the
permutations equal to the identity and to $(12\cdots k)$ simultaneously
minimize the determinant of the minor $M$. This logic is reversible,
proving the equivalence of (1) and (3).

\end{proof}

If the $r$ points of $V$ are in general position, the tropical complex they generate is 
called a 
\emph{generic tropical complex}.
These polyhedral complexes are then polar to the complexes of interior faces of regular
triangulations of $\Delta_{r-1}\times \Delta_{n-1}$.

\begin{cor}
All tropical complexes generated by $r$ points in general position in $\TP^{n-1}$ have the same 
$f$-vector. 
Specifically, the
number of faces of dimension $k$ is equal to the multinomial coefficient

\[
\bin{r+n-k-2}{r-k-1, n-k-1, k} = \frac{(r+n-k-2)!}{(r-k-1)!\cdot (n-k-1)!\cdot k!}.
\]
\end{cor}

\begin{proof}

By Proposition~\ref{thegenericcase}, these objects are in bijection with regular triangulations of
$P=\Delta_{r-1}\times \Delta_{n-1}$.  The polytope $P$ is equidecomposable~\cite{BCS}, meaning that all of
its triangulations have the same $f$-vector. The number of faces of dimension $k$ of the tropical complex
generated by given $r$ points is equal to the number of interior faces of codimension $k$ in the
corresponding triangulation. Since all triangulations of all products of simplices have the same
$f$-vector, they must also have the same interior $f$-vector, which can be computed by taking the
$f$-vector and subtracting off the $f$-vectors of the induced triangulations on the proper faces of $P$.
These proper faces are all products of simplices and hence equidecomposable, so all of these induced
triangulations have $f$-vectors independent of the original triangulation as well.

To compute this number, we therefore need only compute it for one tropical complex. Let the vectors
$v_i$, $1\le i\le r$, be given by $v_i = (i, 2i, \cdots, ni)$. By Theorem~\ref{ineqs}, to count the faces
of dimension $k$ in this tropical complex, we enumerate the existing types with $k$ degrees of
freedom. Consider any index $i$. We claim that for any $x$ in the tropical convex hull of $\{v_i\}$, the
set $\{j\mid i\in S_j\}$ is an interval $I_i$, and that if $i<m$, the intervals $I_m$ and $I_i$ meet in at
most one point, which in that case is the largest element of $I_m$ and the smallest element of $I_i$.

Suppose we have $i\in S_j$ and $m\in S_l$ with $i<m$. Then we have by definition $v_{ij}-x_j \le
v_{il}-x_l$ and $v_{ml}-x_l \le v_{mj}-x_j$. Adding these inequalities yields $v_{ij}+v_{ml} \le
v_{il}+v_{mj}$, or $ij+ml\le il+mj$. Since $i<m$, it follows that we must have $l\le j$. Therefore, we can
never have $i\in S_j$ and $m\in S_l$ with $i<m$ and $j<l$. The claim follows immediately, since the $I_i$
cover $[1,n]$.

The number of degrees of freedom of an interval set $(I_1, \ldots, I_r)$ is easily seen to be the number
of $i$ for which $I_i$ and $I_{i+1}$ are disjoint. Given this, it follows from a simple combinatorial
counting argument that the number of interval sets with $k$ degrees of freedom is the multinomial
coefficient given above. Finally, a representative for every interval set is given by $x_j = x_{j+1}-c_j$,
where if $S_j$ and $S_{j+1}$ have an element $i$ in common (they can have at most one), $c_j=i$, and if
not then $c_j = (\text{min}(S_j)+\text{max}(S_{j+1}))/2$. Therefore, each interval set is in fact a valid
type, and our enumeration is complete. 
\end{proof}

\begin{cor}
The number of combinatorially distinct generic tropical complexes
generated by $r$ points in $\TP^{n-1}$ equals the number of
distinct regular triangulations of $\Delta_{r-1}\times \Delta_{n-1}$.
The number of respective symmetry classes 
under the natural action of the product of symmetric groups
$G=S_r\times S_n$ on both spaces is also the same.
\end{cor}

\begin{figure}
 \begin{center}  \includegraphics{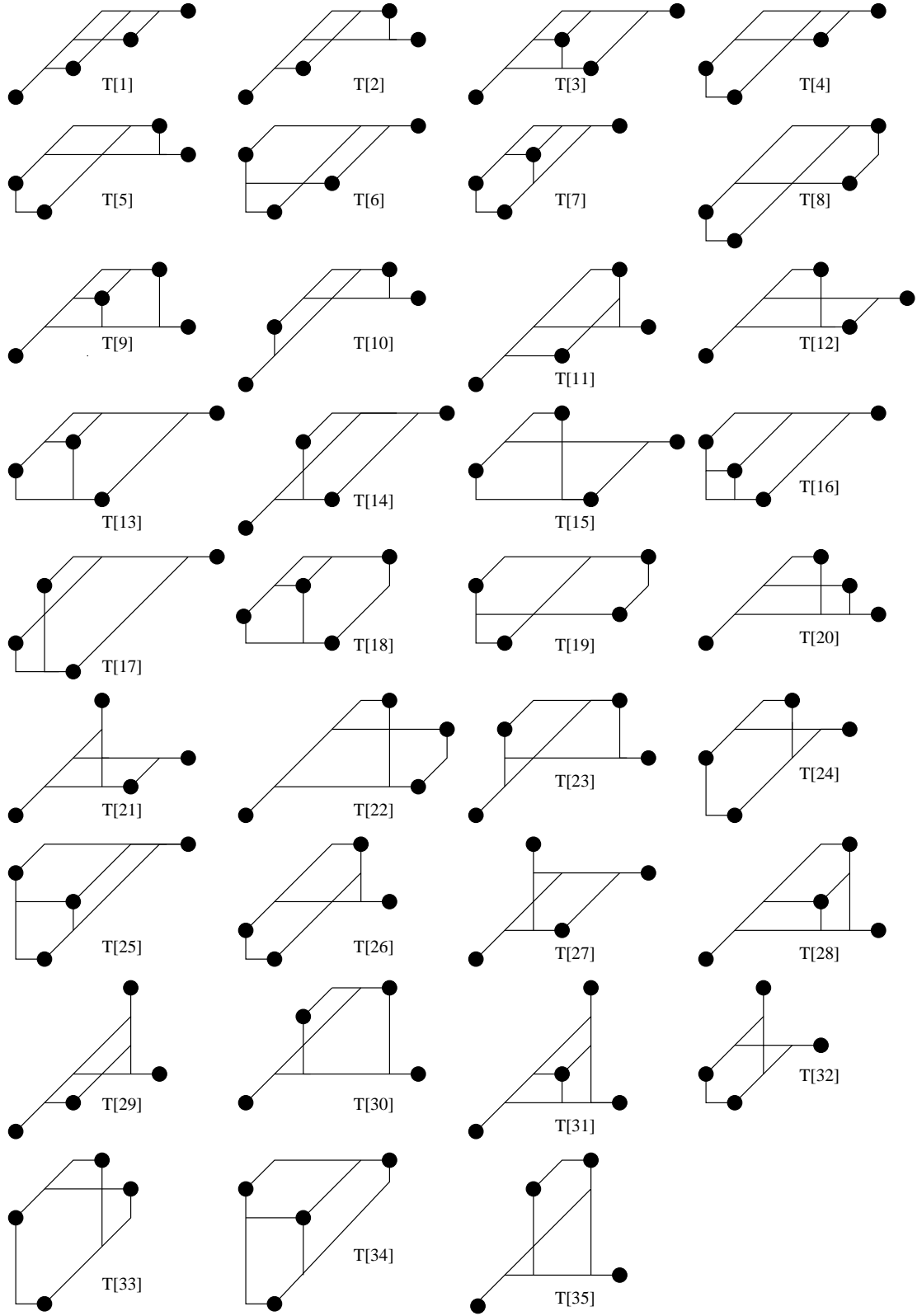}\end{center}
\caption{\label{35quads}
The 35 symmetry classes of tropical complexes generated by four points in $\TP^2$.}
\end{figure}

The symmetries in the group $G$ correspond to a natural action on 
$\Delta_{r-1}\times \Delta_{n-1}$ given by permuting the vertices of the two
component simplices; the symmetries in the symmetric group $S_r$ correspond to permuting the points in a 
tropical polytope, while those in the symmetric group 
$S_n$ correspond to permuting the coordinates. (These are dual, as per Corollary~\ref{polar}.) The number of 
symmetry classes of regular triangulations of the polytope 
$\,\Delta_{r-1} \times \Delta_{n-1}\,$ is 
computable via J\"{o}rg Rambau's TOPCOM~\cite{TOPCOM} for small $r$ and $n$:
$$
\begin{array}{r||r|r}
& 2 & 3 \\
\hline\hline 2 & 5 & 35\\
\hline 3 & 35 & 7,955\\
\hline 4 & 530 &\\
\hline 5 & 13,631\\
\end{array}
$$
For example, the $(2,3)$ entry of the table divulges that there are 
$35$ symmetry classes of regular triangulations of
$\Delta_2\times \Delta_3$. These correspond to the $35$ 
 combinatorial types of four-point configurations in $\TP^2$,
or to the 35 combinatorial types of three-point configurations in $\TP^3$. 
These $35$ configurations (with the tropical complexes they generate) are shown in
Figure~\ref{35quads}; the labeling corresponds to Rambau's 
labeling (see \cite{TOPCOM}) of the regular
triangulations of $\Delta_3\times \Delta_2$.

\section{Phylogenetic analysis using tropical polytopes} \label{trees}

A fundamental problem in bioinformatics is the reconstruction
of phylogenetic trees from approximate distance data. 
In this section we show how tropical convexity
might help provide new algorithmic tools for this problem.
Our approach augments the results in \cite[\S 4]{SS}
and it provides a tropical interpretation of
the work on $T$-theory by Andreas Dress and his 
collaborators \cite{DHM},  \cite{DMT}, \cite{DT}.

Consider a symmetric $n \times n$-matrix $\, D = (d_{ij})\,$
whose entries $\, d_{ij}\,$ are non-negative real numbers
and whose diagonal entries $\, d_{ii} \,$ are all zero.
We say that $D$ is a (finite) \emph{metric} if the
triangle inequality $\, d_{ij} \, \leq \, d_{ik} + d_{jk}\,$ holds 
for all indices $i,j,k$.
Our starting point is the following easy observation:

\begin{prop} \label{metricspace}
The symmetric matrix $D$ is a metric
if and only if all principal $3 \times 3$-minors
of the negated symmetric matrix $\,-D = (-d_{ij})\,$
are tropically singular.
\end{prop}

\begin{proof}
Both properties involve only three points,
so we may assume $n =3$, in which case
$$ - D \quad = \quad 
\begin{pmatrix} 
\phantom{-} 0 &- d_{12} & - d_{13} \\
- d_{12} &\phantom{-} 0 & - d_{23} \\
- d_{13} & - d_{23} &\phantom{-} 0 
\end{pmatrix}. $$
The tropical determinant of this matrix is the minimum
of the six expressions
$$ 0, \, - 2 d_{12}, \, - 2 d_{13}, \, - 2 d_{23}, \, 
-d_{12}-d_{13}-d_{23}\, \, \, \hbox{and} \,\,
-d_{12}-d_{13}-d_{23} . $$
This minimum is attained twice if and only if
it is attained by the last two (identical) expressions,
which occurs if and only if
the three triangle inequalities are satisfied.
\end{proof}

In what follows we assume that $D = (d_{ij})$ is
a metric.
Let $\,P_D \,$ denote the
tropical convex hull in $ \TP^{n-1}\,$
 of the $n$ row vectors (or column vectors)
of the negated  matrix $\, -D = (-d_{ij})$.
Proposition \ref{metricspace} tells us
that the tropical polytope $P_D$ is always one-dimensional 
for $n = 3$.

The finite metric $D = (d_{ij})$ is said to be
a \emph{tree metric}
if there exists a weighted tree $T$ with $n$ leaves
such that
$d_{ij}$ denotes the distance between the
$i$-th leaf and the $j$-th leaf along the
unique path between these leaves in $T$.
The next theorem characterizes tree metrics
among all metrics by the dimension of the
tropical polytope $P_D$. It is the tropical interpretation
of results that are quite classical and well-known in the
phylogenetics literature.

\begin{thm}\label{buneman}
For a given finite metric $D = (d_{ij})$
the following conditions are equivalent:
\begin{enumerate}
\item $D$ is a tree metric,
\item the tropical polytope $P_D$ has dimension one,
\item all $4\times 4$-minors of the matrix $-D$ are tropically singular,
\item all principal $4\times 4$-minors of the matrix $-D$ are tropically 
singular,
\item For any choice of four indices $i,j,k,l \in \{1,2,\ldots,n\}$,
the maximum of the three numbers
$\,d_{ij} + d_{kl},\, d_{ik} + d_{jl}$ and $\, d_{il} + d_{ik}\,$
is attained at least twice.
\end{enumerate}
\end{thm}

\begin{proof}
The condition (5) is the familiar
\emph{Four Point Condition} for tree metrics.
The equivalence of (1) and (5) is a classical
result due to various authors, including
Buneman \cite{Bun} and Zaretsky \cite{Zar}.
See equation (B3) on page 57 in \cite{DHM}.

Suppose that the condition (5) holds.
By the discussion in \cite[\S 4]{SS}, this means that
$\, -D \,$ is a point in the
tropical Grassmannian of lines, in symbols
$\, -D \in {\rm Gr}(2,n) \,\subset \,\TP^{\binom{n}{2}}$.
By \cite[Theorem 3.8]{SS}, the point $-D$
corresponds to a tropical line $L_D $ in 
$\,\TP^{n-1}$. The $n$  distinguished points
whose coordinates are the rows of $-D$ lie on the
line $L_D$. By  Corollary \ref{planesconvex}, it follows
that their tropical convex hull $P_D$ is
contained in $L_D$. This means that
$P_D$ has dimension one, that is, (2) holds.

Suppose that (2) holds. Then the tropical rank of the matrix 
$-D$ is equal to two, by \cite[Theorem 4.2]{DSS}.
This means that all $r \times r$-minors of $-D$
are tropically singular for $r \geq 3$.
The case $r = 4$ is precisely the statement (3).

Obviously, the condition (3) implies
the condition (4). What remains is to prove the
implication from (4) to (5).
For this we note  that the tropical determinant of
the $4 \times 4$-matrix
$$ \begin{pmatrix}
 0 & -d_{12} & - d_{13} & -d_{14} \\
-d_{12} & 0 & - d_{23} & - d_{24} \\
-d_{13} & -d_{23} & 0 & - d_{34} \\
-d_{14} & -d_{24} & -d_{34} & 0 
\end{pmatrix}
$$
equals twice the minimum of 
$\,-d_{12} - d_{34},\, -d_{13} - d_{24}$ and $\, -d_{14} - d_{23}$.
(It's the tropicalization of a $4 \times 4$-Pfaffian).
The matrix is tropically singular
if and only if the minimum is attained twice.
\end{proof}

If the five equivalent conditions of
Theorem \ref{buneman} are satisfied
then the metric tree $T$ coincides with the 
one-dimensional tropical polytope $P_D$. 
To make sense of this statement, we regard
tropical projective space $\TP^{n-1}$ 
as a metric space with respect to the 
infinity norm induced from $\RR^n$,
$$ || x - y || \quad  = \quad
{\rm max} \bigl\{\, | x_i + y_j - x_j - y_i | \,\, :\,\,
1 \leq i < j \leq n \,\bigr\}, $$
and we note that the finite metric 
$D$ embeds isometrically into $P_D$ via the rows of $- \frac{1}{2} D$:
$$ i \,\,\,\, \mapsto \,\,\, \,
\frac{1}{2} \cdot ( -d_{i1}, -d_{i2}, -d_{i3},\ldots, -d_{in}) 
\qquad \hbox{for} \,\, i = 1,2,\ldots,n $$

We learned from  \cite{DMT} that the tropical polytope $P_D$
first appeared in the 1964 paper \cite{Isb} by John Isbell. 
For the proof of the following result
we assume familiarity with results from
 \cite{DHM} and \cite{DMT}.

\begin{thm} The tropical polytope $P_D$ equals 
Isbell's injective hull of the metric $D$. 
\end{thm}

\begin{proof}
According to Lemma \ref{triangtrop},
the tropical polytope $P_D$ is the bounded complex of
the following unbounded polyhedron in  the 
$(2n-1)$-dimensional space $\,W = \RR^{2n}/ \RR(1,\ldots,1,-1,\ldots,-1)$:
$$
{\mathcal P}_{-D} \quad = \quad
\bigl\{\, (y,z) \in W \,:\,\, y_i + z_j \leq - d_{ij} \,\,\,
\hbox{for all} \,\, 1 \leq i,j \leq n \, \bigr\}.
$$
Dress et al.~\cite{DHM} showed that 
the injective hull $T(D)$ of the finite metric $D$
coincides with the complex of bounded faces of
the following $n$-dimensional unbounded polyhedron:
$$
{\mathcal Q}_{-D} \quad = \quad
\bigl\{\, x \in \RR^n \,:\,\, x_i + x_j \geq d_{ij} \,\,\,
\hbox{for all} \,\, 1 \leq i,j \leq n \, \bigr\}.
$$
What we need to show is that the two polyhedra
have the same bounded complex.

The metric $D$ satisfies the tropical
matrix identity $\,-D \, = \, D \odot (-D) $, because
$\, - d_{ij} \, = \, {\rm min}_k (d_{ik} - d_{kj})$.
This implies that any column vector $y$ of $-D$
satisfies $\, y \, = \,  (-y) \odot (-D) $.

Consider any vertex $(y,z)$ of ${\mathcal P}_{-D}$.
Then $y$ is a column vector of $-D$.
Equation (\ref{troplinalg}) implies
$\,z \, = \, (-y) \odot (-D) \,\, = \, \, y$.
Hence  every vertex of ${\mathcal P}_{-D}$
lies in the subspace defined by $y=z$, and 
so does the complex of bounded faces of
${\mathcal P}_{-D}$. Therefore the linear map
$\,(y,z) \mapsto - y \,$ induces
an isomorphism between the bounded
complex of $\,{\mathcal P}_{-D} \,$
and the bounded complex of
$\,{\mathcal Q}_{-D} $.
\end{proof}

Theorem \ref{polar} specifies an involution on the
set of all tropical complexes. We are interested
in the fixed points of this canonical involution.
A necessary condition is that $r=n$ and 
$V$ is a symmetric matrix. 
The previous result and its proof
can be reinterpreted as follows:

\begin{cor} A tropical complex $P$ is
pointwise fixed under the canonical involution
(on the set of all tropical complexes)
if and only if $P$ is the injective hull of
a metric on $\{1,2,\ldots,n\}$.
\end{cor}

\begin{proof}
In order for $P$ to be fixed under the canonical involution,
it is necessary that $n = d$. Hence we can write
$P = {\rm tconv}(-D)$ for some non-negative square matrix $D$.
Now, $P$ is fixed under the involution if and only if
the identity $\,-D \, = \, D \odot (-D) \,$ holds.
This identity is equivalent to $D$ being a metric.
\end{proof}

Dress, Huber and Moulton \cite{DHM} emphasize that
the tropical polytope $P_D$ records many important
invariants of a given finite metric $D$.
For instance, the dimension of $P_D$
gives information about how far the  metric
is from being a tree metric. In practical biological
applications of phylogenetic reconstruction,
the distances $d_{ij}$ are not known exactly,
and $P_D$ appears to contain many ofs the various trees
which are found by existing  software for
phylogenetic reconstruction.

The dimension of the tropical complex $P_D = {\rm tconv}(-D)$
can be characterized combinatorially
by tropicalizing the sub-Pfaffians of
a skew-symmetric $n \times n$-matrix.
The tropical Pfaffians of format $4 \times 4$
specify the four point condition (5) in Theorem \ref{buneman},
while the tropical sub-Pfaffians of format $6 \times 6$
specify the six-point condition which is discussed in
\cite[page 25]{DHM}. The combinatorial study
of \emph{$k$-compatible split systems} can be interpreted
in the setting of tropical algebraic
geometry (cf.~\cite{Mi}, \cite{RGST}, \cite{SS}) as the study of 
the $k$-th \emph{secant variety} in the  
Grassmannian $\, {\rm Gr}(2,n) \subset
\TP^{\binom{n}{2}}$.

Tropical convexity provides a convenient language
to study numerous extensions of the classical
problem of tree reconstruction. As an example,
imagine the following scenario, which would 
correspond to the Grassmannian of planes in $\TP^{n-1}$,
denoted $\, {\rm Gr}(3,n)$.

Suppose  there are
$n$ taxa, labeled $1,2,\ldots,n$, and rather
than having a distance for any pair $i,j$, we
are now given a proximity measure $d_{ijk}$ for any
triple $i,j,k \in \{1,2,\ldots,n\}$. We can then
construct a tropical polytope by taking the
tropical convex hull of $\binom{n}{2}$ points as follows:
$$ P \quad = \quad {\rm tconv}\,\bigl\{ \,
 \bigl(-d_{ij1}, -d_{ij2}, -d_{ij3}, \ldots-d_{ijn} \bigr)
\in \TP^{n-1} \,\,:\,\,
\, 1 \leq i < j \leq n \,\bigr\}. $$
Under certain hypotheses,
the  tropical polytope $P$ can be realized
as the complex of bounded faces of the
polyhedron in $\RR^n$ defined by the inequalities 
$\,x_i + x_j + x_k \geq d_{ijk}$.
It provides a polyhedral
model for the tree-like nature of the 
 data $(d_{ijk})$. The case of 
most interest is when $P$ is two-dimensional
in which case it plays the role of
a \emph{two-dimensional phylogenetic tree}.

The construction of this particular tropical polytope $P$
was pioneered by Dress and Terhalle in the important paper \cite{DT}.
There they discuss \emph{valuated matroids}, which are
essentially the points on the tropical Grassmannian
of \cite{SS}, and they call  $P$ the \emph{tight span of 
a valuated matroid}. We share their view that these
tropical polytopes constitute a promising tool for
 phylogenetic analysis.

\vskip 2cm

\section*{Acknowledgements}
This work was conducted while Mike Develin held
the AIM Postdoctoral Fellowship 2003-2008 and 
Bernd Sturmfels held the MSRI-Hewlett Packard 
Professorship 2003/2004. Sturmfels also acknowledges
partial support from the National Science Foundation
(DMS-0200729).

\end{document}